\documentclass[12pt, oneside]{amsart}
\usepackage{amsfonts}
\usepackage{amsmath,amssymb, amsthm,amsxtra}
\usepackage{amscd} 
\usepackage{verbatim}  %%% comment 

\newcommand{\grassmann}{{\mathbb G_X({d}, \mathcal E)}}

\newcommand{\ronbuntitle}{%
Degree Formula for Grassmann Bundles}
\title[\ronbuntitle]{\ronbuntitle}

\author{Hajime KAJI$^*$ and Tomohide TERASOMA$^{**}$}

\address{%
Department of Mathematics, 
School of Science and Engineering, 
Waseda University 
\newline 
\indent 
3-4-1 Ohkubo, Shinjuku, Tokyo 169--8555, JAPAN }
\email{kaji@waseda.jp}

\address{%
Department of Mathematical Science,
University of Tokyo
\newline 
\indent 
3-8-1 Komaba, Meguro, Tokyo 153-8914, JAPAN}
\email{terasoma@ms.u-tokyo.ac.jp}

\thanks{
\noindent $^*$
Department of Mathematics, 
School of Science and Engineering, 
Waseda University.  
\newline 
\indent 
3-4-1 Ohkubo, Shinjuku, Tokyo 169--8555, JAPAN. 
\newline 
\indent 
{\it E-mail address:} {\tt kaji@waseda.jp}.
\\
\indent 
$^{**}$
Department of Mathematical Science,
University of Tokyo. 
\newline 
\indent 
3-8-1 Komaba, Meguro, Tokyo 153-8914, JAPAN. 
\newline 
\indent 
{\it E-mail address:} {\tt terasoma@ms.u-tokyo.ac.jp}%
\newline
\indent 
2015/04/14}


\theoremstyle{plain}
	\newtheorem{theorem}{Theorem}[section]
	\newtheorem*{theorem*}{Theorem}
	\newtheorem{corollary}[theorem]{Corollary}

\theoremstyle{definition}

\theoremstyle{remark}
	\newtheorem{remark}[theorem]{Remark}

\numberwithin{equation}{section}

\begin{document}

\maketitle

\begin{abstract} 
Let 
$X$ be a non-singular quasi-projective variety over a field,  
and let $\mathcal E$ 
be a vector bundle over $X$. 
Let $\grassmann$ be 
the Grassmann bundle of $\mathcal E$ over $X$ 
parametrizing corank $d$ subbundles of $\mathcal E$
with projection $\pi : \grassmann \to X$, 
let 
$ \mathcal Q \gets \pi^*\mathcal E$ 
be the universal quotient bundle of rank $d$, 
and denote by $\theta$ 
the Pl\"ucker class of $\grassmann$, 
that is, 
the
first Chern class of 
the Pl\"ucker line bundle, 
$\det \mathcal Q$.   
In this short note, a closed formula for 
the push-forward of 
powers of the Pl\"ucker class $\theta$ is 
given in terms of the Schur polynomials in Segre classes of $\mathcal E$, 
which yields a degree formula for 
 $\grassmann$ 
with respect to $\theta$ 
when $X$ is projective and $\wedge ^d \mathcal E$ is very ample. 
\end{abstract}

\setcounter{section}{-1}

\section{Introduction}
 
Let $X$ be a non-singular quasi-projective variety of dimension $n$ 
defined over a field of arbitrary characteristic, 
and let $\mathcal E$ be a vector bundle of rank $r$ over $X$. 
Let $\grassmann$ be 
the Grassmann bundle of $\mathcal E$ over $X$ 
parametrizing corank $d$ subbundles of $\mathcal E$ 
with projection $\pi : \grassmann \to X$,  
and 
let $\mathcal Q \gets \pi^*\mathcal E$ 
be the universal quotient bundle of rank ${d}$ on $\grassmann$.
We denote by $\theta$ 
the first Chern class 
$c_1(\det\mathcal Q)= c_1(\mathcal Q)$
of $\mathcal Q$, and call $\theta$ the 
 {\it Pl\"ucker class} of $\grassmann$. 
Note that the determinant bundle $\det \mathcal Q$ is 
isomorphic to 
the pull-back of 
the tautological line bundle 
$\mathcal O_{\mathbb P_X(\wedge^{d} \mathcal E)}(1)$ 
of $\mathbb P_X(\wedge^{d} \mathcal E)$ by the relative
Pl\"ucker embedding over $X$.

The purpose of this short note is to give 
a closed formula for $\pi_{*}\theta^{N}$,  
the push-forward  of 
powers $\theta^N$ of the Pl\"ucker class $\theta$ to $X$ by $\pi$,  
in terms of the Schur polynomials in Segre classes of $\mathcal E$, 
where $\pi_{*} : A^{*+d( r - d )}(\grassmann)\to A^{*}(X)$ is the push-forward by $\pi$ between 
the Chow rings.

The result is

\begin{theorem}\label{theorem:main_theorem}
For each integer $N \ge d(r-d)$, we have 
$$
\pi_* \theta ^N 
= \sum_{ \vert \lambda \vert = N-d(r-d)} 
f^{\lambda+\varepsilon}
\varDelta _{\lambda}(s(\mathcal E))   
$$
in 
$A^{N-d(r-d)}(X)$, 
where 
$\lambda =(\lambda_1 , \dots, \lambda_d)$ is a partition 
with $\vert \lambda \vert := \sum _{1 \le i \le d} \lambda_{i}$, 
$\varDelta_{\lambda}(s(\mathcal E)):= \det[s_{\lambda_i+j-i}(\mathcal E)]_{1 \le i,j\le d}$ is  
the Schur polynomial in Segre classes of $\mathcal E$ corresponding to  
$\lambda$, 
$\varepsilon := (r-d)^d = (r-d, \dots, r-d)$, 
and 
$f^{\lambda + \varepsilon}$ is 
the number of standard Young tableaux with shape $\lambda+\varepsilon$. 
\end{theorem}

The Segre classes $s_i(\mathcal E)$ here are
the ones satisfying
$s(\mathcal E)c(\mathcal E^{\vee})=1$ as in 
\cite{fujita}, \cite{laksov}, \cite{laksov-thorup}, 
where 
$s(\mathcal E)$ and $c(\mathcal E)$ denote respectively 
the total Segre class and the total Chern class of $\mathcal E$. 
Note that our Segre class $s_i(\mathcal E)$ 
differs by the sign $(-1)^i$ 
from the one in \cite{fulton}.

\begin{corollary}[degree formula for Grassmann bundles]%
\label{corollary:degree_formula}
If $X$ is projective 
and 
$\wedge^{d} \mathcal E$ is very ample, then  $\grassmann$  
is embedded in the projective space $\mathbb P(H^0(X, \wedge^{d} \mathcal E))$ 
by the tautological line bundle $\mathcal O_{\grassmann}(1)$,  
and its degree is given by 
$$
\deg \grassmann = 
 \sum_{ \vert \lambda \vert = n} 
f^{\lambda+\varepsilon}
\int_X \varDelta _{\lambda}(s(\mathcal E)) 
. 
$$
\end{corollary}

Here a vector bundle $\mathcal F$ over $X$ is said to be {\it very ample} if 
the tautological line bundle $\mathcal O_{\mathbb P_X(\mathcal F)}(1)$ 
of $\mathbb P_X(\mathcal F)$ is very ample.

Setting  $n:=0$, we recover 
the degree formula of Grassmann varieties, as follows:

\begin{corollary}[{\cite[Example 14.7.11 (iii)]{fulton}}]
\label{corollary:classical_degree_formula}
Let $\mathbb G({d},  r)$ be the Grassmann variety 
parametrizing codimension  $d$ subspaces %${d}$-dimensional quotient spaces 
of a vector space of dimension $r$.  
Then its degree with respect to the Pl\"ucker embedding is given by
$$
\deg 
\mathbb G({d},  r)
=\frac{{({d}(r -{d}))!} \prod_{1 \le l \le {d}-1} l ! }{\prod_{1 \le l \le {d}} (r - l )! }
.
$$
\end{corollary}

\section{Proofs}

\begin{proof}[{\it Proof of Theorem \ref{theorem:main_theorem}}]
Let $\xi_1 , \dots, \xi_d$ be the Chern roots of the universal quotient bundle $\mathcal Q$. 
Then we can write $\theta = \xi_1 + \cdots + \xi_d$ formally. 
Using Pieri's formula \cite[\S2.2]{fulton-young} repeatedly, 
and applying Jacobi-Trudi identity \cite[I, (3.4)]{macdonald}, 
we obtain that 
$$
\theta^N 
=
\sum_{\vert \mu \vert = N} 
f^{\mu} 
\varDelta_{\mu} (\underline{\xi}) 
=
\sum_{\vert \mu \vert = N} 
f^{\mu} 
\varDelta_{\mu} (s(\mathcal Q)) , 
$$
where 
$\varDelta_{\mu}(\underline{\xi})$
is the Schur polynomial in 
$\underline{\xi}:= (\xi_1, \dots , \xi_d)$ 
corresponding to a partition $\mu$. 
It follows from the push-forward formula of J\'ozefiak-Lascoux-Pragacz
 \cite[Proposition 1]{jlp} 
that 
$$
\pi_*\varDelta_{\mu} (s(\mathcal Q))
= \varDelta_{\mu-\varepsilon} (s(\mathcal E))
. 
$$
Therefore 
we obtain 
$$
\pi_* \theta^N 
= \sum_{\vert \mu \vert = N} 
f^{\mu} 
\varDelta_{\mu-\varepsilon} (s(\mathcal E))
= \sum_{\vert \lambda \vert = N-d(r-d)} 
f^{\lambda+\varepsilon} 
\varDelta_{\lambda} (s(\mathcal E))  , 
$$
where 
$\lambda$ is a partition, 
and $\varepsilon := (r-d)^{d} = (r-d,\dots,r-d)$.  
\end{proof}

\begin{proof}[{\it Proof of Corollary \ref{corollary:degree_formula}}]
By the assumption  
$\grassmann$  
is projective and 
the tautological line bundle $\mathcal O_{\mathbb P_X(\wedge^{d} \mathcal E)}(1)$ 
defines 
an embedding 
$\mathbb P_X(\wedge^{d} \mathcal E) \hookrightarrow 
\mathbb P(H^0(X, \wedge^{d} \mathcal E))$. 
Therefore 
$\grassmann$
is considered to be a projective variety 
in 
$\mathbb P(H^0(X, \wedge^{d} \mathcal E))$ 
via 
the relative Pl\"ucker embedding 
$\grassmann \hookrightarrow \mathbb P_X(\wedge^{d} \mathcal E)$
over $X$ 
defined by the quotient 
$\wedge^d\pi^*\mathcal E \to \wedge^d \mathcal Q=\det \mathcal Q$. 
Since the hyperplane section class of $\grassmann$ is equal to the Pl\"ucker class $\theta$, 
we obtain the conclusion, 
taking $N:=\dim \grassmann = {d} (r -{d}) + {n}$ in Theorem \ref{theorem:main_theorem}. 
\end{proof}

\begin{proof}[{\it Proof of Corollary \ref{corollary:classical_degree_formula}}]
The conclusion follows from Corollary \ref{corollary:degree_formula} with $n:=0$, 
since the number $f^{\lambda+\varepsilon}$ is known to be given as follows 
(\cite[p.53]{fulton-young} and 
\cite[p.54, Exercise 9]{fulton-young}):
\begin{equation*}
f^{\lambda+\varepsilon}  
=
\frac{N! \prod_{1 \le i < j \le d } (\lambda _i - \lambda_j - i+ j ) }{\prod_{1 \le i \le d} (r + \lambda_i - i )!}  
. 
\qedhere 
\end{equation*}
\end{proof}

\begin{remark}
Under the same assumption as in Theorem \ref{theorem:main_theorem},
one can prove a push-forward formula of the following form: 
\begin{equation*}
\pi_* \theta^N 
= 
\sum_{\vert k\vert = N-d(r-d)} 
\frac{N! \prod_{1 \le i<j\le d} (k_i  -k_j -i + j) }
{\prod_{1 \le i \le d}(r + k_i   -i)}
\prod_{i=1}^{d} s_{k_i}(\mathcal E)
\end{equation*}
in 
$A^{N-d(r-d)}(X) \otimes \mathbb Q$, 
where 
$k = (k_1, \dots , k_d) \in \mathbb Z_{\ge 0}^d$  with $\vert k \vert := \sum_i k_i$, 
and $s_i(\mathcal E)$ is the $i$-th Segre class of $\mathcal E$. 

\end{remark}

\medskip  

\noindent{\it Acknowlegements.} 
The authors would like to 
thank the first referee 
for his/her detailed comments and invaluable advice: In fact, 
our original proof 
is much longer than 
and completely different from 
the one given here, which is 
due to the referee. 
The authors thank Professor Hiroshi Naruse and 
Professor Takeshi Ikeda, too, for useful discussion and kind advice. 

The first author is supported by JSPS KAKENHI Grant Number 25400053. 
The second author is supported by JSPS KAKENHI Grant Number 15H02048.%

\end{document}